\documentclass[%
  onecolumn
]{mpi2015-cscpreprint}

\usepackage[american]{babel}

\usepackage{amsmath,amsthm,amssymb}
\usepackage{graphics,graphicx,epsfig,algorithm,algpseudocode}
\usepackage{graphicx}
\usepackage{acronym}
\usepackage{bookmark}
\usepackage{subfigure}

\def\ba{\ensuremath{\mathbf{a}}}
\def\bu{\ensuremath{\mathbf{u}}}
\def\bx{\ensuremath{\mathbf{x}}}
\def\bf{\ensuremath{\mathbf{f}}}
\def\bd{\ensuremath{\mathbf{d}}}
\def\bb{\ensuremath{\mathbf{b}}}
\def\bv{\ensuremath{\mathbf{v}}}
\def\bvcnn{\ensuremath{\mathbf{v}_{\mathsf{CNN}}}}

\def\fc{\ensuremath{{\mathsf{fc}}}}

\def\bw{\ensuremath{\mathbf{w}}}

\def\bp{\ensuremath{\mathbf{p}}}

\def\bz{\ensuremath{\mathbf{z}}}
\def\bbeta{\ensuremath{\boldsymbol{\beta}}}
\def\brho{\ensuremath{\boldsymbol{\rho}}}

\def\bU{\ensuremath{\mathbf{U}}}
\def\bF{\ensuremath{\mathbf{F}}}
\def\bG{\ensuremath{\mathbf{G}}}
\def\bW{\ensuremath{\mathbf{W}}}
\def\bV{\ensuremath{\mathbf{V}}}
\def\bD{\ensuremath{\mathbf{D}}}
\def\bA{\ensuremath{\mathbf{A}}}
\def\hbA{\ensuremath{\hat{\mathbf{A}}}}
\def\bJ{\ensuremath{\mathbf{J}}}
\def\bI{\ensuremath{\mathbf{I}}}
\def\bM{\ensuremath{\mathbf{M}}}
\def\hbM{\ensuremath{\hat{\mathbf{M}}}}
\def\bN{\ensuremath{\mathbf{N}}}
\def\hbN{\ensuremath{\hat{\mathbf{N}}}}
\def\bS{\ensuremath{\mathbf{S}}}

\newacro{POD}[POD]{\emph{Proper Orthogonal Decomposition}}
\newacro{FEM}[FEM]{\emph{Finite Element Method}}
\newacro{AE}[AE]{\emph{autoencoder}}
\newacro{cPOD}[cPOD]{\emph{clustered POD}}
\begin{document}
  
%%%%%%%%%%%%%%%%%%%%%%%%%%%%%%%%%%%%%%%%%%%%%%%%%%%%%%%%%%%%%%%%%%%%%%%%%%%%%%%%
% PAPER INFORMATION.                                                           %
%%%%%%%%%%%%%%%%%%%%%%%%%%%%%%%%%%%%%%%%%%%%%%%%%%%%%%%%%%%%%%%%%%%%%%%%%%%%%%%%

\title{Convolutional Autoencoders, Clustering, and POD for Low-dimensional Parametrization of Navier-Stokes Equations}
  
\author[$\ast$,$\dagger$]{Yongho Kim}
\affil[$\ast$]{Max Planck Institute for Dynamics of Complex Technical Systems, Magdeburg, Germany}%\authorcr

\author[$\ast$,$\dagger$]{Jan Heiland}
\affil[$\dagger$]{Department of Mathematics, Otto von Guericke University Magdeburg,
    Magdeburg, Germany\authorcr
    \email{ykim@mpi-magdeburg.mpg.de}, \orcid{0000-0003-4181-7968}
  \email{heiland@mpi-magdeburg.mpg.de}, \orcid{0000-0003-0228-8522}}
  
\shorttitle{CAEs, Clustering and POD for Low-dimensional Parametrization}
\shortauthor{Y. Kim, J. Heiland}
\shortdate{}
  

\msc{MSC1, MSC2, MSC3}
  
\abstract{
Simulations of large-scale dynamical systems require expensive computations. 
Low-dimensional parametrization of high-dimensional states such as Proper
Orthogonal Decomposition (POD) can be a solution to lessen the burdens by
providing a certain compromise between accuracy and model complexity. 
However, for really low-dimensional parametrizations (for example for
controller design) linear methods like the POD come to their natural limits so 
that nonlinear approaches will be the methods of choice.
% In particular, Linear Parameter Varying (LPV) systems enable to control high-dimensional and nonlinear systems efficiently 
% and the performance of the LPV approach is affected by low-dimensional parametrization. 
% However, linear projection methods such as POD would not guarantee a certain parametrization accuracy in very low-dimensional parametrization cases. 
In this work we propose a convolutional autoencoder (CAE) consisting of a
nonlinear encoder and an affine linear decoder and consider combinations with
k-means clustering for improved encoding performance. 
The proposed set of methods is compared to the standard POD approach in two cylinder-wake scenarios modeled by the incompressible Navier-Stokes equations. 
}

 \novelty{
 \begin{itemize}
\item Convolutional autoencoders (CAEs) that handle FEM data for low-dimensional parametrization.
\item Combination of encoders and decoders with clustering techniques for
  improved performance.
\item Thorough demonstration of the capabilities of CAEs with clustering in
  comparison with POD in two cylinder-wake simulation scenarios.
\end{itemize}
 }

\maketitle

%%%%%%%%%%%%%%%%%%%%%%%%%%%%%%%%%%%%%%%%%%%%%%%%%%%%%%%%%%%%%%%%%%%%%%%%%%%%%%%%
% PAPER CONTENT.                                                               %
%%%%%%%%%%%%%%%%%%%%%%%%%%%%%%%%%%%%%%%%%%%%%%%%%%%%%%%%%%%%%%%%%%%%%%%%%%%%%%%%
  
\section{Introduction}%\label{sec:intro}

Accurate \ac{FEM} discretizations of fluid flow easily
lead to dynamical systems with millions of degrees of freedom. This poses a
computational challenge to simulations and, due to the involved nonlinearities,
even more to computer-aided controller design based on these models.
% With numerical analysis, control systems have been developed successfully in diverse fields such as engineering and chemistry.
% Usually, the accuracy of relates to how many elements are dealt with. 
% However, fluid flow is composed of a huge amount of molecules so it requires high computational costs to produce the flow in a simulation.
% For instance, a liter of water contains more than $10^{25}$ molecules. 
% Therefore, a big challenge is to handle high-dimensional systems efficiently. 

The core idea and the promise of model order reduction techniques is the identification and exploitation of
lower-dimensional coordinates that can well represent or approximate the
dynamics of a generally high-dimensional systems. 
For example, \emph{projection methods} like the successful method of \ac{POD} base on linear
projections onto subspaces that are designed to encode the most relevant
dynamics of the system. 
These methods can efficiently reduce the dimension of the state space of a given
dynamical systems and lead to significant savings of memory and computational
requirements in simulations. However, linear methods are naturally limited in
their accuracy for a given number of degrees of freedom as expressed by the
\emph{Kolmogorov n-width}, see e.g., \cite{OhlR16}.

Accordingly, if a very low-dimensional parametrization of still satisfactory accuracy is
wanted, one may need to resort to nonlinear approaches;
cp. \cite{HeiBB22,LeeC20,SyFaJa18}.

Thus, in view of controller design where a low dimension is weighted higher than the
actual accuracy, we consider methods to approximate and parametrize the states on a
low-dimensional manifold by possibly non-linear maps.

Parametrization refers to representing the states of a given system in 
different coordinates, approximation means that this representation is designed to
trade in accuracy for maximal low-dimensionality of the parametrizing
coordinates. 
By \emph{encoding} we will refer to the (computation of) reduced and
parametrized representation of a state, whereas \emph{decoding} or
\emph{reconstruction} describes the
computation of the state in the actual coordinates from the parametrized
representation. The according mechanisms or algorithms are referred to as
\emph{encoder} or \emph{decoder}.

The method of \ac{POD}, thus, can be seen as an encoder/decoder based on linear
projections. As mentioned above, linear methods are limited
% base of a linear encoder and a linear decoder so that it is applicable to 
% diverse fields as well as it is easy to analyze the linear relation between 
% full-order states and their reduced coordinates. (see, e.g., \cite{Ba18POD,Mi21PODGal})
% However, as POD has unveiled the limitation of the linear mechanism in terms of 
% the reconstruction performance, 
so that 
nonlinear methods such as general \ac{AE}
methods have been emerging as an alternative or an enhancement to POD
\cite{FreM22}. 

An \emph{autoencoder} (see e.g., \cite[Ch. 14]{Ian16DL}) often refers to
a \emph{neural network} that is designed to efficiently encode data and that
is trained by just considering the data itself. 
% an unsupervised learning method that seeks to generate an output that is equal to each input.
The underlying network architecture can be problem specific and contain several
hidden layers. 
Nonlinearity in the encoding and decoding possibly arises from using activation
function.
Based on these principles, many types of autoencoders such as sparse autoencoders, denoising autoencoders, 
contractive autoencoders, and convolutional autoencoders have been developed. (see e.g., \cite{Dor20AE})

With the increasing availability of computing facilities and toolboxes for
design and training of neural networks, since well over a decade, autoencoders
have been considered as alternatives in the field of model order reduction.
A recent example considers a modified nonlinear autoencoder that approximates
the elements of the coefficient matrices of a linearized two-link planar robot
manipulator; see \cite{PaRo20}.
Other examples target the solution of convection-dominated problems through deep
convolutional autoencoders; see, e.g., \cite{LeeC20})

% Above all, in terms of dimension reduction, it is a key factor that a model consists of 
% an encoder extracting low-dimensional features and a decoder reconstructing the input, 
% and when it comes to \emph{linear parameter varying} (LPV) systems, autoencoders are often compared with 
% traditional and powerful dimension reduction methods such as POD
% and \emph{Principal Component Analysis} (PCA).

With PDEs as the underlying model, the variants of autoencoders relying on
so-called \emph{Convolutional Neural Networks} (CNNs) \cite{Lec98CNN}
seem to be particularly useful as indicated by several attempts to use CNNs 
to solve PDEs; see e.g., \cite{Han21phy, KimCh22FCNN, Ni19ConvPDE} for current
works.
The realization of the linear part of the propagation from one layer to another as a convolution
drastically reduces the parameters of the network by enforcing \emph{sparse
connectivity}, meaning that over the layers every variable is related to only a
limited number of other variables. 
Moreover, the accompanying operation of
\emph{pooling}, that coalesces several neighboring variables into one,
reduces the variable dimension in every layer.
In classical \ac{FEM} approaches, sparse connectivity is a major performance
criterion and generally achieved by considering basis functions with local
support whereas the \emph{pooling} can be related to multiscale approaches.

% as well as sparse connectivity makes a model small compared to \emph{Multi-Layer Perceptrons} (MLPs) 
% when they handle the same number of feature maps. 
% have led the computer vision field to successful prosperity. 
% The convolution is basically a \emph{linear transformation} so it can be decided 
% whether a layer is linear or nonlinear depending on the type of activation functions as well as 
% the composite of convolutional layers preserves the linearity of the operation. 

% In this paper, we consider \emph{Convolutional Autoencoders} (CAEs) \cite{Mas11CAE} to improve 
% the reconstruction performance and to exclude POD from the model order reduction procedure.

Once a very low-dimensional code of a state has been derived, it is the task of
the decoder to reconstruct the high-dimensional state from the code. 
Intuitively spoken, the design of the decoder becomes increasingly complex the
higher variations in the output space it should reproduce. 
Accordingly, we suppose that clustering in the reduced order coordinates will
well separate the state space so that
separate individual decoders then only need to cover a certain range in the
overall variety; see \cite{HeiKim22} for an example implementation. 

\emph{Clustering} algorithms are used to distinguish data and divide
a given dataset into certain subsets without the use of labels, i.e. without
knowing characteristics of the subset in advance; see, e.g., \cite{AS21cluster}
for clustering in general and \cite{Scu10Kmeans} for a discussion of the most popular clustering
algorithm -- the \emph{$k$-means clustering} -- that iteratively identifies a
predifined number of clusters in a data set based on distance measures.
Moreover, clustering algorithms well integrate with neural networks, the more
that clustering becomes easier in low data dimensions like in low-dimensional
feature spaces; see e.g., \cite{Ma20Deepk,Ju16clustering,Xi17DeepClu}.

% groups a dataset into $k$ clusters based on the distance
% between each vector and iteratively updated optimally located centroids of the clusters. 

Since clustering is not a smooth operation, such an approach has not been attempted in dynamical systems until recently.
Nonetheless, it has been observed that dimensionality reduction techniques with clustering improve retrieval performance 
and reduce computational costs 
by individually processing subsets of a large-scale dataset, e.g., in text
retrieval systems; see \cite{JG05CSVD}. %  for general considerations.

In this work, with the focus on identifying very low-dimensional parametrizations of
states (as opposed to low-dimensional dynamical models) the nonsmoothness is not
an issue so that we can directly adopt clustering in the numerical studies.

% As an unsupervised learning approach, clustering algorithms can be out of the
% potential risk of unmatched annotations and reduce the considerable time for
% data labeling. 
% Thus, it would be more useful than manual labeling for classifying data that have vague classification criteria such as flows.
% It is not only
% applicable easily with neural networks, but also it has fewer hyperparameters
% than other clustering methods.
% Therefore, we investigate how CAEs and $k$-means clustering can be used for low-dimensional parametrization with FEM data.

This paper is organized as follows:
In \Cref{sec:nseq}, we introduce the incompressible Navier-Stokes equation which will be
our FEM model under consideration.
In \Cref{sec:ldp}, we define low-dimensional parametrization methods.
In \Cref{sec:mor}, we discuss how to build each reduced order model using considered methods respectively.
In \Cref{sec:results}, we show the results of the reconstruction errors of the methods and evaluate the reduced order models based on their residuals.
In \Cref{concl}, we sum up the paper, discuss limitations, and lay out potential
research directions for future work.

%%%%%%%%%%%%%%%%%%%%%%%%%%%%%%%%%%%%%%%%%%%%%%%%%%%%%%%%%%%%%%%%%%%%%%%%%%%%%%%%
\section{Incompressible Navier-Stokes Equations}\label{sec:nseq}
We consider the dynamical system of spatially discretized incompressible Navier-Stokes equations
\begin{subequations}\label{ns}
\begin{align}
\bM\dot{\bv}(t) + \bN(\bv(t))\bv(t)+ \bA\bv(t) - \bJ^\top\bp(t) &= \bf(t) \\
\bJ\bv(t)&=0,
\end{align}
\end{subequations}
where for time $t>0$, $\bv(t)\in\mathbb R^{n_v}$ and $\bp(t)\in \mathbb R^{n_p}$ denote the
states of the velocity and the pressure respectively on the FEM
mesh, and where $\bM\in\mathbb R^{n_v\times n_v}$ is the mass matrix, where $\bA\in\mathbb R^{n_v\times n_v}$ models the
diffusion, where $\bN\colon \mathbb R^{n_v}\to \mathbb R^{n_v\times n_v}$
represents the discretized convection, and where $\bJ\in \mathbb R^{n_p\times
n_v}$ and
$\bJ^\top$ denote the discrete divergence and gradient, respectively; see
\cite{BehBH17} for technical details and example discretizations.

In what follows, we will drop the time dependency in the variables $\bv$ (and
later also $\brho$ and $\tilde \bv$ which will denote the code and the reconstruction).

As it is common practice at least for theoretical considerations, the system
\eqref{ns} can be expressed as an ordinary differential equation. Here, we
will use the ODE formulation to avoid technical difficulties with treating the
$p$ variable in numerical schemes (cp. \cite{AltH15}) and to be directly
adaptable to other ODE models.

We briefly lay out, how the ODE formulation of the incompressible Navier-Stokes
equations is obtained. Under the reasonable assumption that $\bJ$ is of full rank and
$\bM$ is invertible and symmetric positive definite, by means of  the equation
$\bp=\bS^{-1}\bJ\bM^{-1}(\bN(\bv)\bv + \bA\bv-\bf)$ where
$\bS=\bJ\bM^{-1}\bJ^\top$, we can eliminate 
$\bp$ from \eqref{ns} and obtain
\begin{equation*}%\label{ns1}
\bM\dot{\bv} + (\bI-\bJ^\top\bS^{-1}\bJ\bM^{-1})(\bN(\bv)\bv + \bA\bv - \bf)=0.
\end{equation*}
With the projector $\Pi:=\bI-\bM^{-1}\bJ^\top\bS^{-1}\bJ$ and having confirmed that 
the solution $v$ satisfies $v=\Pi v$, we can state that for the computation of
$v$, system \eqref{ns} can be replaced by
\begin{equation}\label{ns2}
\bM\dot{\bv} + \Pi^\top(\bN(\bv)\bv + \bA\bv - \bf)=0.
\end{equation}
% By the \emph{state-dependent coefficient} (SDC) matrix $\bN(\bv)$, the nonlinear convection term of the equation can be computed efficiently in the form of $\bN(\bv)\bv$.

In the next sections we will investigate how well the state $\bv$ can be encoded
in a very low-dimensional code $\rho$ with respect to several measures. 
Firstly, the standard reconstruction error
\[
  \|\bv(t)-\tilde \bv(t)\|_{\bM}
\]
that we
measure in the so-called $\bM$-norm that is induced by the mass matrix $\bM$ and
that is discrete and consistent counterpart of the $L^2$ \emph{Sobolev}
norm associated with the underlying model problem. 
Secondly, the reconstruction with respect to the resulting
convection is considered by measuring
\[
  \|\bN(\tilde \bv(t))\bv(t) - \bN(\bv(t))\bv(t)\|_{\bM^{-1}},
\]
where the $\bM^{-1}$ norm is chosen to account for the FEM context in which $\bN(\bv(t))\bv(t)$ is a functional.

Finally, in view of possibly providing related surrogate low-order dynamical
models, we will report on the residuals at the data points.

% However, the equation (\ref{ns2}) requires a high computational cost to handle the flows as a \emph{full order model} (FOM). \emph{Reduced order model} (ROM) can lower the computational requirement while guaranteeing a certain accuracy. Moreover, $\bN(\bv)$ can be implemented by a few parameters $\rho_i$, $i=1,2,\cdots ,n_\rho$ (e.g., $n_\rho=2,3$) as a \emph{linear parameter varying} (LPV) system if $\bN(\bv)$ is linear and an approximation $\bN(\bv) \approx \tilde{\bN}(\brho(\bv))$ is obtained with a low-dimensional parameter $\brho$.

%%%%%%%%%%%%%%%%%%%%%%%%%%%%%%%%%%%%%%%%%%%%%%%%%%%%%%%%%%%%%%%%%%%%%%%%%%%%%%%%
\section{Low-dimensional Parametrization}\label{sec:ldp}
In this section, we consider a number of possibly nonlinear methods for encoding
and decoding with a focus on very low-dimensional parametrizations and with the
standard method POD as a benchmark.

As introduced above, we will use $\bV$ to denote the velocity with values in
$\mathbb R^{n_v}$, $\brho$ to denote the code with values in $\mathbb R^{n_\rho}$,
and $\tilde \bv$ to denote the reconstruction of the velocity from the code.
The considered approaches are as follows:

\subsection{Proper Orthogonal Decomposition (POD)}\label{subsec:pod}
The method of \ac{POD} \cite{Ber93POD} is readily interpreted as an autoencoder with 

\begin{itemize}
  \item[(1)] a
\emph{linear} encoder defined as
\begin{equation*}
\brho=\bV^\top\bv 
\end{equation*}
\item[(2)]
and a \emph{linear} decoder 
\begin{equation*}
\tilde{\bv}=\bV\brho,
\end{equation*}
\end{itemize}
where $\bV\in\mathbb{R}^{n_v \times n_\rho}$ is a POD basis that is, typically,
computed as the $n_\rho$ leading singular vectors of a matrix of snapshots of
the variable $\bv$.

% , $\bv\in\mathbb{R}^{n_v}$ is a velocity state obtained by a FEM simulation, and $\brho\in\mathbb{R}^{n_\rho}$ is a low-dimensional parameter $(n_v \gg n_\rho)$.

\subsection{Convolutional Neural Network (CNN)}\label{subsec:cnn}
In order to apply standard convolutional approaches, the input data needs to be
available on a tensorized grid, i.e., a grid of rectangles with no hanging nodes
and that are aligned with the coordinate axes; see \cite[Fig. 3 and 4]{HeiKim22}
for an illustration of this interpolation.
To achieve that make use of corresponding interpolation matrix $\bI_c$ such that
the CNN input $\bvcnn(t)$ can be generated by the linear transformation
$\bvcnn(t)=\bI_c\bv(t)$. We consider the CNN architecture introduced as introduced in
\cite{HeiBB22} that makes use of POD modes for the reconstruction. Thus, 
the CNN is defined as a model consisting of
\begin{itemize}
  \item[(1)] 
 a \emph{nonlinear convolutional} encoder $\mu$:
\begin{equation*}
\brho=\mu(\bvcnn)
\end{equation*}
\item[(2)] 
and a \emph{linear} decoder
\begin{equation*}
\tilde{\bv}=\phi(\brho)=\bV(\bW\brho),
\end{equation*}
\end{itemize}
where $\bvcnn\in\mathbb{R}^{2\times w \times h}$ is a CNN input, $\bV\in\mathbb{R}^{n_v \times r}$ is a POD basis with $r$ modes, and $\bW\in\mathbb{R}^{r \times n_\rho}$ is a trainable matrix. 
The feedforward network of the decoder $\mu$ in (1) is concretely described as follows:
\begin{subequations}\label{cnnnet}
\begin{align*}
&\bu^{(m)} = \ba(\bF^{(m)}(\bu^{(m-1)})),\, \bu^{(0)}=\bvcnn,\quad m=1,\cdots, M-1,\\
&\bu^{(M)} = \bF^{(M)}(\bu^{(M-1)}),\\
&\brho=\bU\bar\bu^{(M)}+\bbeta,\quad  \bU\in\mathbb{R}^{n_\rho\times n_u},\, \bbeta\in\mathbb{R}^{n_\rho},
\end{align*}
\end{subequations}
where $\ba$ is a nonlinear activation function, $M$ is the number of convolutional layers, $\bF$ is a convolutional layer, $\bar\bu^{(M)}$ is the vectorized $\bu^{(M)}$, and $n_u$ is the dimension of $\bar\bu^{(M)}$. 

%That is, $\mu(\bvcnn)= U\bar F^{(M)}(a(F^{(M-1)}(\cdots a(F^1(\bvcnn)))+\bbeta$

\subsection{Convolutional Autoencoder (CAE)}\label{subsec:cae}

\begin{figure}[t]
\centering
\includegraphics[width=1.0\columnwidth]{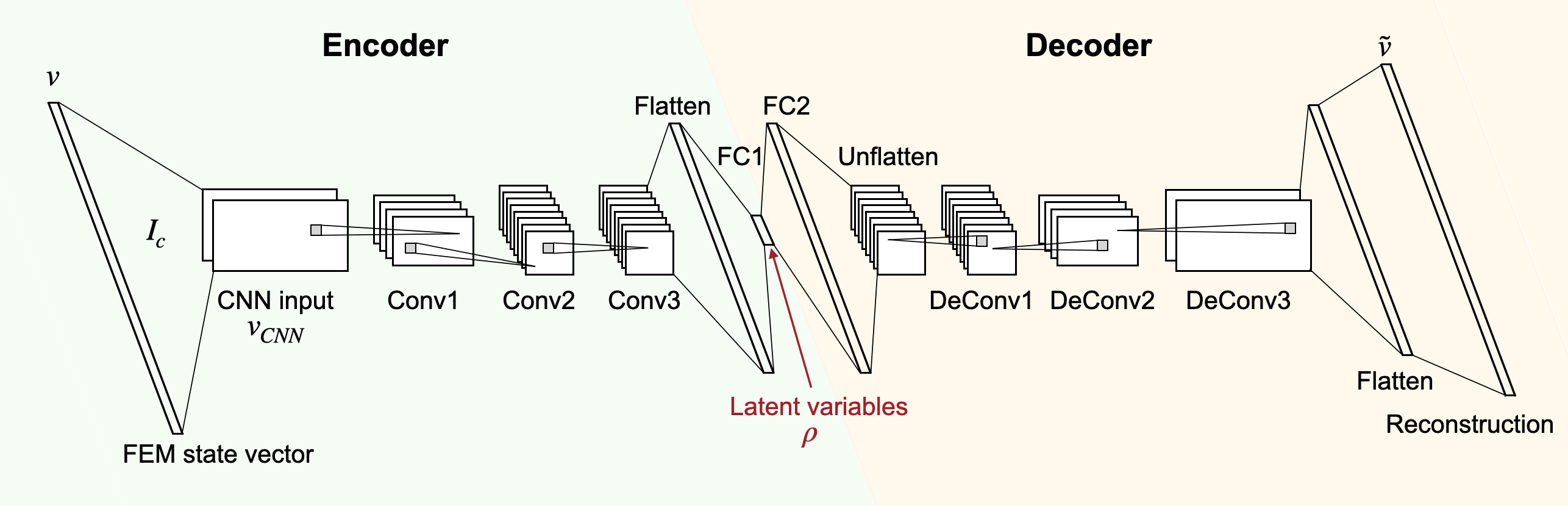}
\caption{Convolutional Autoencoder for the reconstruction of FEM state vectors}%
\label{fig:cae-architecture}
\end{figure}

% Data constructed by rectangular grids are basically required to feed inputs into CNNs. Our previous work \cite{HeiKim22} demonstrated that CAEs outperform POD using CNN data interpolated by $I_c$. In this paper, we train CAEs using FEM data directly. 
Other than using POD modes for the reconstruction, we propose the use of
so-called \emph{transposed} convolutional layers without nonlinear activation
and with an interpolation operator $\bI_p$ that interpolates the values from the
tensorized grid back to the FEM grid. Thus, we can investigate the CNN
properties for reconstruction while -- as a composition of
affine-linear operators -- keeping the process affine-linear.

This proposed model, thus, consists of

\begin{itemize}
\item[(1)] a \emph{nonlinear convolutional} encoder
\begin{equation*}
\mu \colon \bvcnn \to \brho
\end{equation*}
and
  \item[(2)] an \emph{affine linear deconvolutional} decoder 
\begin{equation*}
\tilde{\bv}=\phi(\brho)=\bD\brho+\bb
\end{equation*}
where $\bD\in\mathbb{R}^{n_v\times n_\rho}$ is a matrix, and $\bb\in\mathbb{R}^{n_v}$ is a vector. 
\end{itemize}
In practice, the decoder $\phi \colon \brho \to \tilde \bv$ is realized via:
\begin{subequations}
\begin{align*}
% &\bu^{(m)} = \ba(\bF^{(m)}(\bu^{(m-1)})),\, \bu^{(0)}=\bvcnn, m=1,\cdots, M,\\
% &\brho=\bU_1\bar\bu^{(M)}+\bbeta_1,\,  U_1\in\mathbb{R}^{n_\rho\times n_u}, \bbeta_1\in\mathbb{R}^{n_\rho},\\
&\bar\bz^{(0)}=\bU\brho+\bbeta,\quad  U\in\mathbb{R}^{n_u\times n_\rho}, \bbeta\in\mathbb{R}^{n_u},\\
&\bz^{(m)} = \bG^{(m)}(\bz^{(m-1)}),\, m=1,\cdots, M, \quad \text{with $\bz^{(0)}$ as the unflattened $\bar\bz^{(0)}$}\\
&\tilde{\bv} = \bI_p\bar\bz^{(M)} ,\quad  \bI_p\in\mathbb{R}^{n_z\times n_v},
\end{align*}
\end{subequations}
where $M$ is the number of layers, $\bG$ is a deconvolutional layer, $\bar\bz^{(M)}$ is the vectorized $\bz^{(M)}$, and $n_z$ is the dimension of $\bar\bz^{(M)}$.
The architecture is shown in \Cref{fig:cae-architecture}. 

\subsection{Clustered POD (cPOD)}\label{subsec:clustered-pod}
As motivated in the introduction, the reconstruction performance can well be improved by using multiple decoders depending on clusters. 
% Moreover, it is expected that their centroids and labels can be used to build controller selectors in a feedback control system. 
Before we go for clustering with general autoencoders, we define a
\ac{cPOD} as a model consisting of
\begin{itemize}
  \item[(1)] a \emph{linear} encoder
\begin{equation*}
\brho=\bV^\top\bv,
\end{equation*}
\item[(2a)] a cluster selection algorithm
  \begin{equation}\label{eq:cluster}
    c\colon \brho(t) \mapsto l \in \{1,2,\dotsc,k\},
  \end{equation}
\item[(2b)] and \emph{$k$ linear} decoders 
\begin{equation*}
\tilde{\bv}=\bV_{l}(\bV_{l}^\top\bV)\brho,
\end{equation*}
\end{itemize}
where $\bV$ is a POD basis based on the whole data, where
$\bV_{l}\in\mathbb{R}^{n_v \times n_\rho}$ is a POD basis associated with the
$l$-th cluster, $l=1,2,\cdots,k$ and where $k$ is the number of clusters. 
For setting up the model, we employ $k$-means clustering on the low-dimensional
state vectors $\brho$, classify the actual data $\bv$ accordingly and then
compute POD bases for each cluster separately. 
Thus, during evaluation, a label $l \in \{1, \dotsc, k\}$ is extracted from
the reduced state vector $\brho(t)$ red and a proper \ac{POD} basis is selected
for the given label. 
Note that clustering is \emph{nonlinear} and \emph{discontinuous} so the
decoding becomes nonlinear and discontinuous as well.

\subsection{Individual CAE (iCAE)}\label{subsec:icae}
As another clustering approach, we adapt the individual CAE approach (iCAEs)
that we proposed in \cite{HeiKim22} for CNNs with POD based reconstruction for
the use of general CAEs as in \Cref{subsec:cae}. An iCAE model consists of
\begin{itemize}
\item[(1)] a \emph{nonlinear convolutional} encoder
\begin{equation*}
\mu\colon \bvcnn \rightarrow \brho
\end{equation*}
as in \Cref{subsec:cae}(1), of (2a) a clustering operation as in
\Cref{subsec:clustered-pod}(2a) and
 \item[(2b)]
\emph{$k$ affine linear deconvolutional} decoders 
\begin{equation*}
\tilde{\bv}=\phi_{l}(\brho),
\end{equation*}
\end{itemize}
where $\phi_l$ is the $l$-th affine linear deconvolutional decoder trained by the $l$-th cluster. The encoder is continuous and differentiable. The decoding, on the other hand, is nonlinear and discontinuous. It is known that data compression can improve clustering accuracy. Since a CAE is trained in advance of training an iCAE, the pretrained encoder can be used without any extra training, and reduced low-dimensional states $\brho$ can be used for $k$-means clustering. Thus, we train only $k$ decoders based on $k$ clusters while freezing the encoder parameters. In the inference, a reduced vector $\brho$ and a label are yielded simultaneously once the encoder extracts a reduced state vector, and then a proper decoder is selected based on the label.

%%%%%%%%%%%%%%%%%%%%%%%%%%%%%%%%%%%%%%%%%%%%%%%%%%%%%%%%%%%%%%%%%%%%%%%%%%%%%%%%
\section{Low-order Parametrization of the State and the Models}\label{sec:mor}
\begin{table}[t]
 \begin{center}
  \caption{Description of the CNN-based models}\label{cnndesc}
   \begin{tabular}{c c c}
   \hline
       $-$ & CNN & CAE \\ 
     \hline 
     $\#$encoding layers & $4+\fc$ & $3+\fc$ \\
     $\#$convolution channels & $(4,8,10,12)$ & $(4,8,8)$ \\
     $\#$decoding layers & $\fc+\bV$ &  $\fc+3+\fc$ \\
     $\#$deconvolution channels & - & $(8,8,4)$ \\
     \hline
     \multicolumn{3}{l}{\footnotesize * $\fc$: a fully connected layer, $\bV$: a POD basis with $r$ modes}\\
   \end{tabular}
 \end{center}
\end{table} 

\begin{figure}[t]
\centering
\includegraphics[width=1\columnwidth]{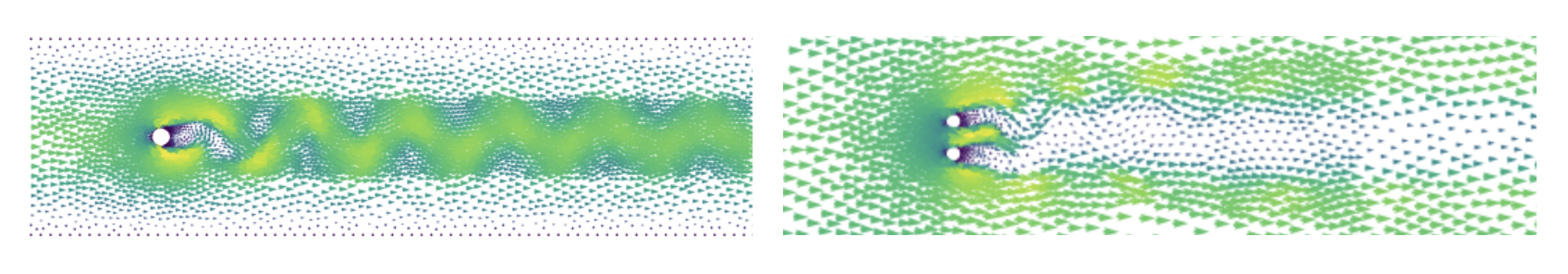}
\caption{Two snapshots of developed velocity states: (left) a single cylinder case at Re=40 (right) a double cylinder case at Re=50}%
\label{fig:devsnap}
\end{figure}

\begin{figure}[t]
\centering
\includegraphics[width=1\columnwidth]{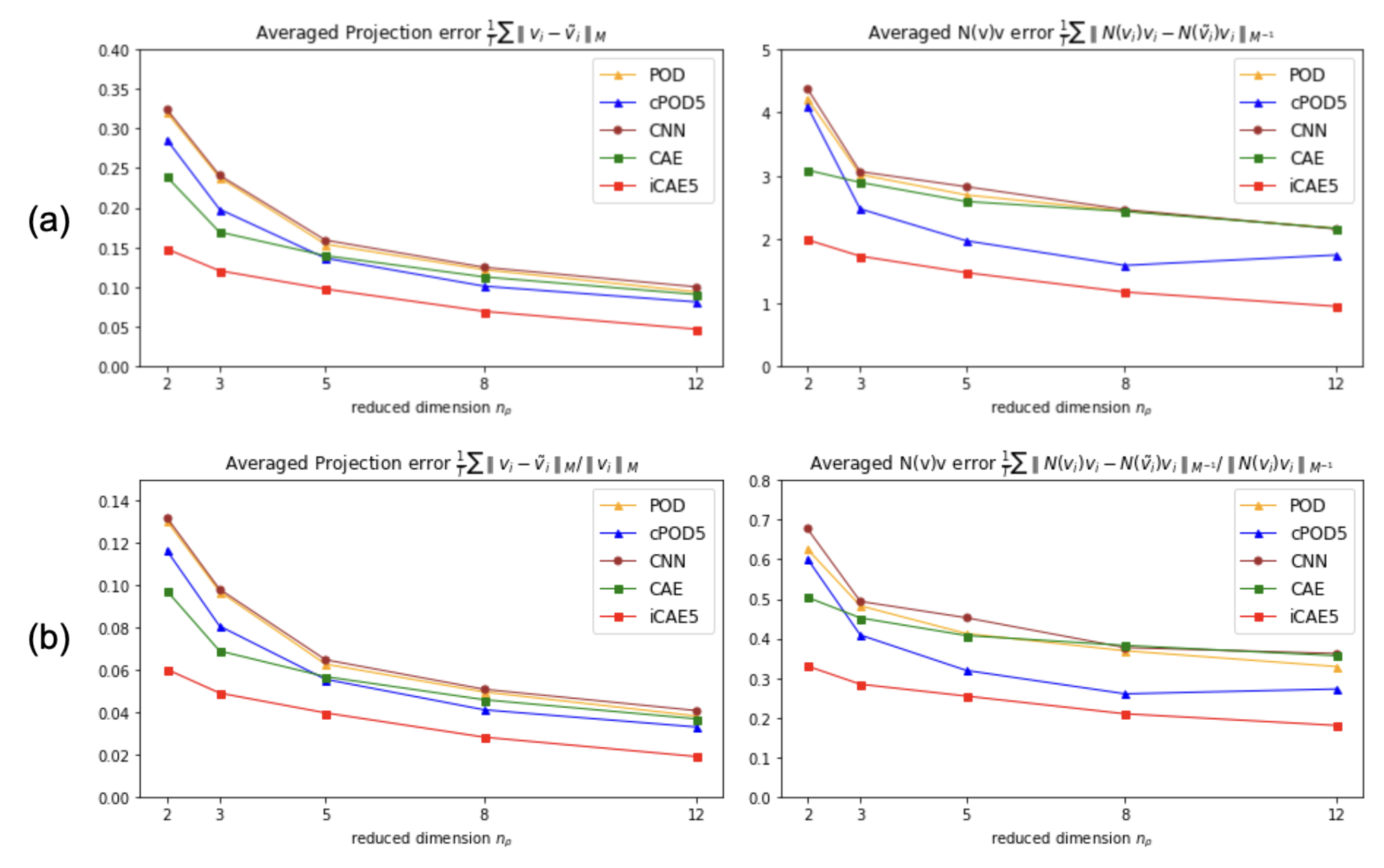}
\caption{(a) Averaged reconstruction errors (b) Averaged relative errors for the single cylinder case (\Cref{sec:num-setup-single-cyl}).}%
\label{fig:rec_err-single}
\end{figure}

In this section, we apply the aforementioned models of \Cref{sec:ldp} to define low-dimensional
LPV approximations of the incompressible Navier-Stokes equation in the ODE form
of \eqref{ns2}.
Where applicable, i.e., where the decoding is linear, we lay out the
affine-linear structure of the LPV approximation that is naturally induced by
the quadratic nature of the Navier-Stokes equations since $\bv \mapsto \bN(\bv)$
is linear; see also \cite{HeiBB22}.
% As reduced order models, we apply the aforementioned models and affine Linear
% Parameter Varying (LPV) approximations of $\bN(\bv)$ to the full order model
% (\ref{ns2}) respectively. 
We note that we are mostly interested in low-dimensional parameterized
approximations of the nonlinear term. % of the state equations that we present here. 
Nonetheless, for the methods with a differentiable decoder, we also state the
equations in the reduced order coordinates.
% rather than in a classical reduced order
% model. 
% Therefore we provide 
% \begin{itemize}
%   \item the general LPV approximation, 
%   \item the affine-linear realization if it exists, and
%   \item the model equations in reduced order coordinates.
% \end{itemize}
In particular, as mentioned, the clustering methods lead to non-continuous
reconstructions $\tilde{\bv}$ so that $\dot{\tilde \bv}$ might not be defined. 

% and so
% that it cannot be used to parametrize the time derivative of $\bv$. In this paper,
% We only parametrize the nonlinear convection $\bN(\bv)\bv$ 
% through the methods introduced in Section 
Numerical evaluations of the residuals for two test cases will be reported
% and evaluate each residual based on the left side of each equation 
in \Cref{sec:results}.

\begin{enumerate}
  \item For the \textbf{POD} encoder, the parametrization reads
\begin{equation*}%\label{eq:nspod}
  \bM\dot{\bv} + \Pi^\top[\bigl (\sum_{i=1}^{n_\rho}\rho_i(\bv)\bN(\bw_i)\bigr)\bv+ \bA\bv - \bf]=0
\end{equation*}
where $\brho(\bv)=\bV^\top\bv$, and $\bw_i$ is the $i$-th column vector of the
POD basis $\bV$ as defined in \Cref{subsec:pod}. The model in the
reduced coordinates reads
\begin{equation}\label{eq:nspod-rom}
  \hbM\dot{\brho} + \Pi^\top[\bigl
(\sum_{i=1}^{n_\rho}\rho_i\hbN(\bw_i)\bigr)\brho+ \hbA\brho - \bf]=0,
\end{equation}
where $\hbM:=\bM\bV$, $\hbA:=\bA\bV$, and $\hbN(\cdot):=\bN(\cdot)\bV$. Further
note that the standard POD reduced order model is obtained from
\eqref{eq:nspod-rom} by a multiplication by $\bV^\top$ from the left.

\item For the \textbf{CNN} method, the parametrized equations read
\begin{equation*}%\label{eq:nscnn}
  \bM\dot{\bv} + \Pi^\top[\bigl(\sum_{i=1}^{n_\rho}\rho_i(\bv)\bN(\bw_i)\bigr)\bv+ \bA\bv - \bf]=0
\end{equation*}
where $\brho(\bv)=\mu(\bv)$, and $\bw_i$ is the $i$-th column vector of
$\bV\bW$; see \Cref{subsec:cnn}. The equations in the reduced coordinates are
obtained by
\begin{equation*}%\label{eq:nscnn}
  \hbM\dot{\brho} +
  \Pi^\top[\bigl(\sum_{i=1}^{n_\rho}\rho_i\hbN(\bw_i)\bigr)\brho+ \hbA\brho - \bf]=0
\end{equation*}
where $\hbM:=\bM\bV\bW$, $\hbA:=\bA\bV\bW$, and $\hbN(\cdot):=\bN(\cdot)\bV\bW$.

\item For the \textbf{CAE} model, the parametrized equations read
\begin{equation*}%\label{eq:nscae}
\bM\dot{\bv} + \Pi^\top[(\bN(\bb)+\sum_{i=1}^{n_\rho}\rho_i(\bv)\bN(\bd_i))\bv+ \bA\bv - \bf]=0
\end{equation*}
where $\brho(\bv)=\mu(\bv)$ and $\bb$ is the bias of the decoder and $\bd_i$ is the $i$-th column vector of $\bD$; see Section \ref{subsec:cae}.
The equations in the reduced coordinates are
obtained by
\begin{equation*}
\hat{\bM}\dot{\brho} + \Pi^\top[(\bN(\bb)+\sum_{i=1}^{n_\rho}\rho_i\bN(\bd_i))(\bD \brho + \bb)+ \hat{\bA}\brho + \hat{\bf}]=0
\end{equation*}
where $\hbM:=\bM\bD$, $\hbA:=\bA\bD$, and $\hat{\bf}:=\bA\bb - \bf$.

\item For the \textbf{cPOD} model, the parametrized equations read
\begin{equation*}%\label{eq:nscpod}
\bM\dot{\bv} + \Pi^\top[(\sum_{i=1}^{n_\rho}\rho_i(\bv)\bN(\bw_i^l))\bv+ \bA\bv - \bf]=0
\end{equation*}
where $\brho(\bv)=\bV^\top\bv$ and $\bw_i^{l}$ is the $i$-th column vector of the
matrix $\bV_{l}(\bV_{l}^\top\bV)$ with the overall POD basis $\bV$ and the POD
basis $\bV_l$ of the corresponding cluster; see Section \ref{subsec:clustered-pod}.

\item For the \textbf{iCAE} model, the parametrized equations read
\begin{equation*}%\label{eq:nsicae}
\bM\dot{\bv}  + \Pi^\top[(\bN(\bb_l)+\sum_{i=1}^{n_\rho}\rho_i(\bv)\bN(\bd_i^{l}))\bv + A\bv - \bf]=0
\end{equation*}
where $\brho(\bv)=\mu(\bv)$, $\bb_l$ is the bias of the $l$-th decoder and $\bd_i^l$ is the $i$-th column vector of $\bD_l$; see Section \ref{subsec:icae}.
\end{enumerate}

As for the \textbf{iCAE} and \textbf{cPOD} models, we note that the dependence
of $l=l(\brho)$ on $\brho$ via the labelling algorithm \eqref{eq:cluster}
introduces a nonlinearity in the
decoding. Nonetheless, as the $n_\rho \cdot k$ coefficients $\bN(\bw_i^l)$ or
$\bN(\bd_i^l)$ can be precomputed, $i=1, \dots, n_\rho$, $l=1, \dotsc, k$, where $k$ is the number
of clusters, the inherent affine linear structure in the reduced order models
can be exploited for efficient realizations in numerical simulations.
% $l$ depends on $\brho$ so that the nonlinear dependency can be described as $l=l(\brho)$.
% For convenience, we use the notation $l$ for the clusters.  so that the nonlinear dependency can be described as $l=l(\brho)$.
% For convenience, we use the notation $l$ for the clusters. 

%%%%%%%%%%%%%%%%%%%%%%%%%%%%%%%%%%%%%%%%%%%%%%%%%%%%%%%%%%%%%%%%%%%%%%%%%%%%%%%%
\section{Numerical Examples}\label{sec:results}
\begin{figure}[t]
\centering
\includegraphics[width=1\columnwidth]{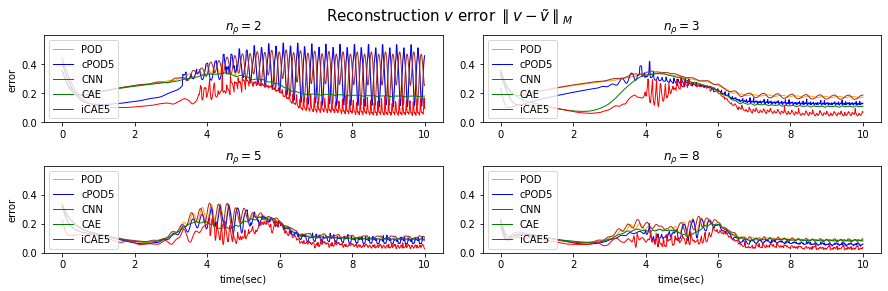}
\caption{Reconstruction error in $[0,10]$ for the single cylinder case (\Cref{sec:num-setup-single-cyl}).}%
\label{fig:v_err-single}
\end{figure}

\begin{figure}[t]
\centering
\includegraphics[width=1\columnwidth]{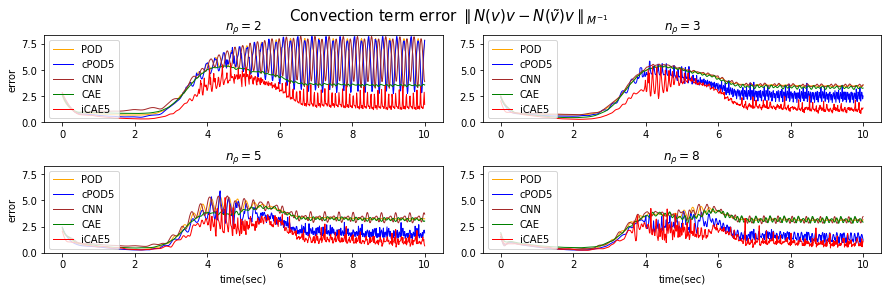}
\caption{Convection error in $[0,10]$ for the single cylinder case (\Cref{sec:num-setup-single-cyl}).}%
\label{fig:nvv_err-single}
\end{figure}

\begin{figure}[t]
\centering
\includegraphics[width=1\columnwidth]{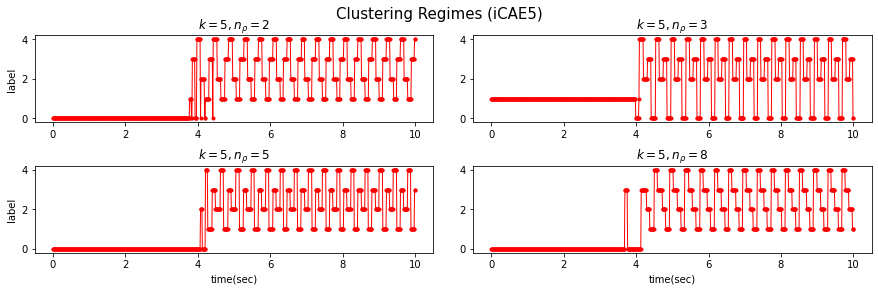}
\caption{Prompt clustering reaction in $[0,10]$ for the single cylinder case (\Cref{sec:num-setup-single-cyl}).}%
\label{fig:clt-single}
\end{figure}

\begin{figure}[t]
\centering
\includegraphics[width=1\columnwidth]{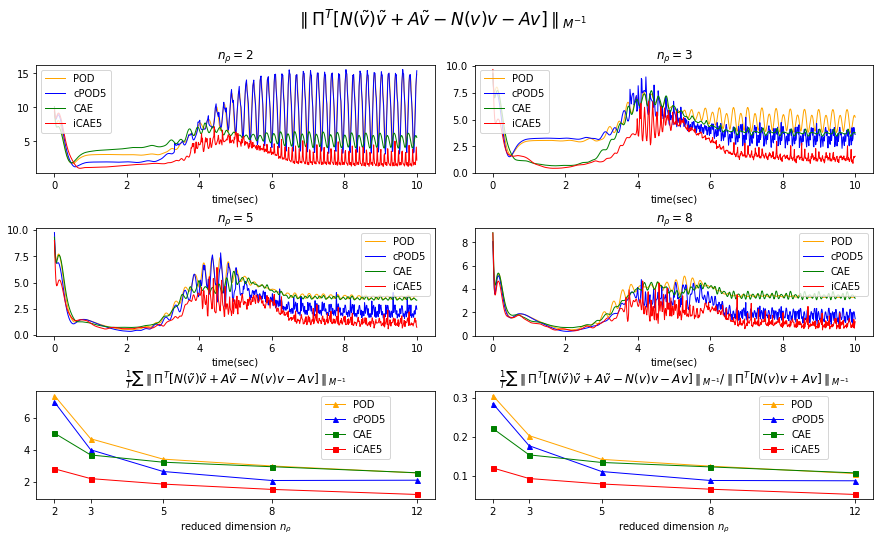}
\caption{Error of $ \Pi^\top[\bN(\tilde{\bv})\tilde{\bv}+ \bA\tilde{\bv} - (\bN(\bv)\bv+ \bA\bv )]$ in $[0,10]$ (\Cref{sec:num-setup-single-cyl}).}
\label{fig:residual-single}
\end{figure}

The cylinder wake is an established benchmark example that generates diverse
flow patterns depending on simulation settings and the Reynolds number. 
In this section, we consider two-dimensional cylinder-wake phenomena observed
behind a single cylinder and two cylinders respectively as shown in \Cref{fig:devsnap}. 
We investigate how the models generate reconstructed state
vectors by using data including periodic wakes behind a single cylinder and data
containing chaotic vortices behind two cylinders respectively. 
Moreover, the main topic of the work is to investigate very low-dimensional
parametrization of incompressible flows so $n_\rho=2,3,5,8,12$ are used.

The raw data, all routines and a script for reproducing the presented results
are publicly available as described at the end of the manuscript in 
\Cref{sec:code-availability}.

\newpage  % to keep this sentence together with the list that follows
For the evaluation of all experiments, we use the six evaluation metrics as
follows:
\begin{itemize}
\item The averaged reconstruction error
 \begin{equation*}%\label{vmetric}
% E(\bv,\tilde{\bv})=
 \frac{1}{T}\sum_{i=1}^{T}\parallel \bv_i-\tilde{\bv}_i\parallel_\bM,
 \end{equation*}
 \item the averaged relative error of reconstruction $\bv$
 \begin{equation*}%\label{vmetric-relative}
%  E(\bv,\tilde{\bv})=
 \frac{1}{T}\sum_{i=1}^{T}\parallel \bv_i-\tilde{\bv}_i\parallel_\bM/\parallel \bv_i\parallel_\bM
 \end{equation*}
\item the averaged convection term error,
\begin{equation*}%\label{nvvmetric}
% E(\bv,\tilde{\bv})=
\frac{1}{T}\sum_{i=1}^{T}\parallel
\bN(\bv_i)\bv_i-\bN(\tilde{\bv}_i)\bv_i\parallel_{\bM^{-1}},
\end{equation*}
\item the averaged relative error of reconstruction $\bN(\bv)\bv$
\begin{equation*}%\label{nvvmetric-relative}
% E(\bv,\tilde{\bv})=
\frac{1}{T}\sum_{i=1}^{T}\parallel
\bN(\bv_i)\bv_i-\bN(\tilde{\bv}_i)\bv_i\parallel_{\bM^{-1}}/\parallel
\bN(\bv_i)\bv_i\parallel_{\bM^{-1}},
\end{equation*}
\item the averaged residual % $\Pi^\top(\bN(\bv)\bv+ \bA\bv-\bf)$
\begin{equation*}%\label{metric:residual}
% R(\bv,\tilde{\bv})=
\frac{1}{T}\sum_{i=1}^{T}\parallel \Pi^\top[(\bN(\tilde{\bv})\tilde{\bv}+ \bA\tilde{\bv}) - (\bN(\bv)\bv+ \bA\bv )]\parallel_{\bM^{-1}},
%\frac{1}{T-1}\sum_{i=2}^{T}\parallel \bM \delta^-\bv_i +
%\Pi^\top(\bN(\tilde{\bv}_i)\bv_i + \bA\bv_i - \bf) \parallel_{\bM^{-1}},
\end{equation*}
\item and the averaged relative residual %of $\Pi^\top(\bN(\bv)\bv+ \bA\bv-\bf)$
\begin{equation*}%\label{metric:relresidual}
% R(\bv,\tilde{\bv})=
\frac{1}{T}\sum_{i=1}^{T}\parallel \Pi^\top[(\bN(\tilde{\bv})\tilde{\bv}+ \bA\tilde{\bv}) - (\bN(\bv)\bv+ \bA\bv )]\parallel_{\bM^{-1}}/\parallel \Pi^\top(\bN(\bv)\bv+ \bA\bv)\parallel_{\bM^{-1}},
%\frac{1}{T-1}\sum_{i=2}^{T}\parallel \bM \delta^-\bv_i +
%\Pi^\top(\bN(\tilde{\bv}_i)\bv_i + \bA\bv_i - \bf) \parallel_{\bM^{-1}},
\end{equation*}
\end{itemize}
where $\parallel\bx\parallel_\bM=\sqrt{\bx^\top\bM\bx}$, $\bv_i$ is a target
velocity, $\tilde{\bv}_i$ is a reconstructed velocity and $T$ is the number of
snapshots.
%, $\bN(\cdot)$ is the state-dependent coefficient matrix, and $\delta^-
%{\bv}_i$ denotes the backward difference quotient $\frac{\bv_i-\bv_{i-1}}{\delta
%t}$ with $\delta_t$ as the time difference between the data points. 

\subsection{Single Cylinder}\label{sec:num-setup-single-cyl}
\begin{figure}[t]
\centering
\includegraphics[width=1\columnwidth]{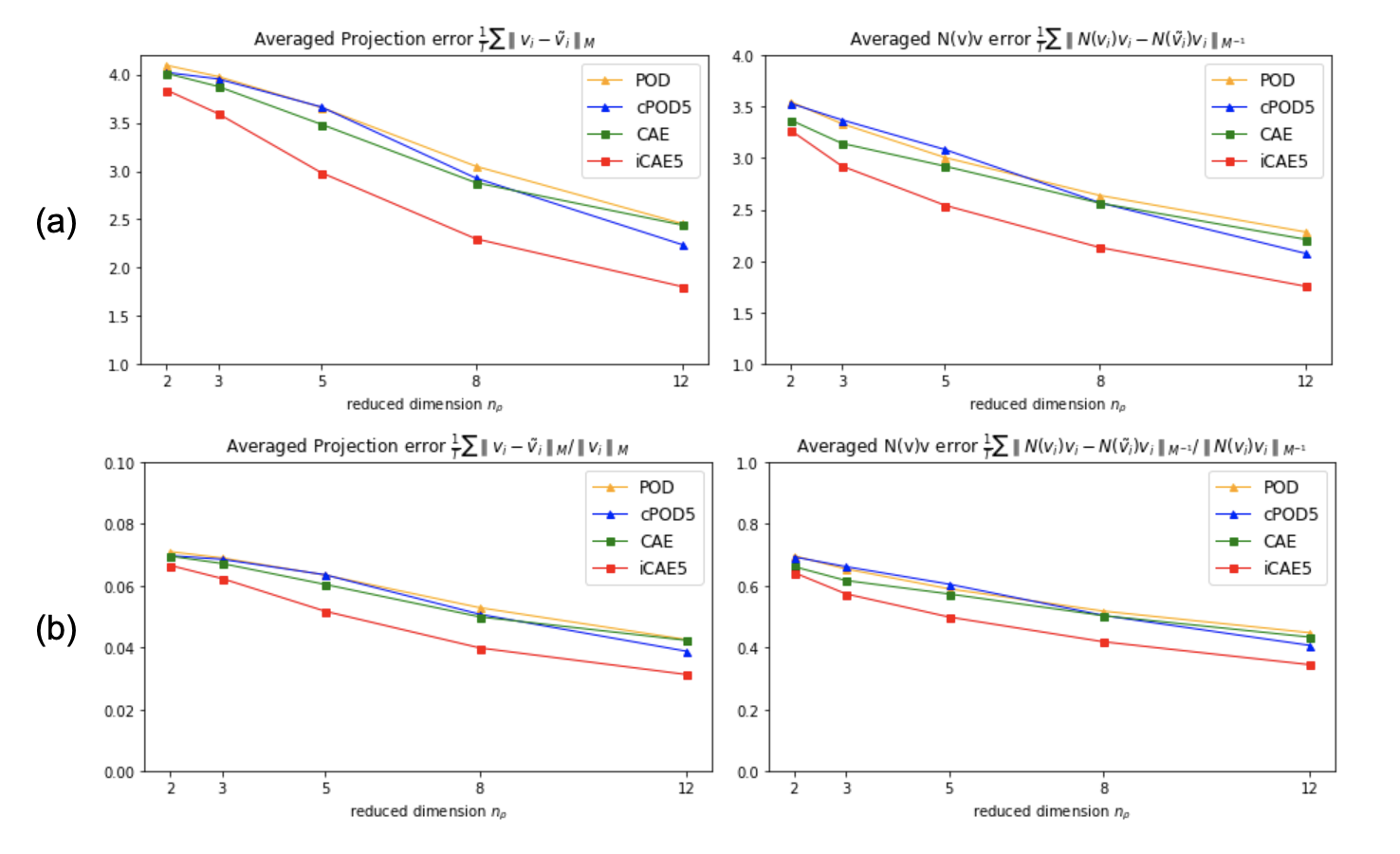}
\caption{(a) Averaged reconstruction errors and (b) Averaged relative errors for
the double cylinder case (\Cref{sec:num-setup-double-cyl}).}
\label{fig:rec_err-double}
\end{figure}

As for the training, namely to set up the POD basis and to train the neural
networks and the clustering algorithm, we use 400 equidistant snapshots of $x$ and $y$ directional velocities in
$[0,10]$. The snapshots were generated by a FEM simulation of Equation
(\ref{ns}) for a cylinder wake with a Reynolds number of 40 starting at the
somewhat nonphysical associated Stokes solution.
All data points $\bv(t)\in\mathbb{R}^{42764}$ (i.e., $n_v=42764$) are spatially
distributed in the domain $(0,5)\times (0,1)$. To feed the data into the
convolutional neural networks, the CNN input data $\bvcnn(t) \in \mathbb
R^{2\times 47\times 63}$
% ($400\times2\times47\times63$) 
are generated by means of an interpolation matrix
$\bI_c\in\mathbb{R}^{42764\times 5922}$ followed by a reshaping of the data
vector into the 3D tensor. 
For the evaluation dataset, the same conditions were used to generate 800 snapshots on $[0,10]$.

The parameters of the CNN/CAE models are listed in \Cref{cnndesc}. 
Basically, the CNN model includes four convolutional layers and a fully
connected layer in the encoder, and the decoder uses a linear transformation and
a POD basis $\bV\in\mathbb{R}^{42764\times15}$ (i.e., $r=15$) to map the reduced
parameters back to the original state space. 

The encoder of the CAE contains three convolutional layers and a fully connected
layer. The decoder conducts successively a fully connected layer, three
deconvolutional layers, and a fully connected layer. In the encoder part, a nonlinear activation function \emph{ELU} \cite{Djo16ELU} is applied in each hidden layer of the encoders. For the clustering models, the number of clusters $k$ is 5 and the considered clustering models are called cPOD5 and iCAE5 respectively. 

In the training session of the networks, we use the ADAM optimizer
\cite{Kin15ADAM} with a multistep learning rate scheduler and a batch size of 64. 
As a loss function, the mean square error is used. When the CAE is trained using
$\bvcnn$ data, an interpolation matrix $\bI_p$ satisfying $\bv=\bI_p\bvcnn$ is
used to return to the states from the tensorized grid to the FEM mesh.
Note that because of the rather rough tensorized grid, the two interpolations
mean a major loss of details and that, in particular, $\bv \neq \bI_p \bI_c \bv$.
% , due to the different coarseness is no $\bI_p$ which is the inverse of $\bI_c$ so $\bv \approx \bI_p\bvcnn$ should be considered. 

As the use of $\bI_p$ and $\bI_c$ means a significant loss of accuracy, 
one could try to rather learn these interpolations operations.
Hence, we optimize $\bI_p$ in the last layer of the decoder and other decoder parameters simultaneously,
instead of the usage of the fixed $\bI_p$ obtained by the optimization of $\bv \approx \bI_p\bvcnn$.

\Cref{fig:rec_err-single} shows how the cPOD and the iCAE can improve the reconstruction performance of the standard POD and the CAE respectively. 
The CNN model and POD have similar results. In very low-dimensional cases ($n_\rho=2,3$), the CAE outperforms the cPOD. Regarding the convection term errors, a similar trend to the projection errors is observed except for the cPOD. The error graph of the cPOD increases at $n_\rho=12$, contrasting with trends that decrease by getting larger $n_\rho$. \Cref{fig:v_err-single} shows that the POD and the cPOD suffer from the large performance degradation when $n_\rho=2$. Taken overall, the iCAE is the best over time. However, peaks occur at around 4 seconds when $n_\rho=5$.

\Cref{fig:nvv_err-single} indicates that the $\bN(\bv)\bv$ errors tend to
be similar to the reconstruction $\bv$ errors and that the pretrained models can
approximate $\bN(\bv)\bv$ equally well, although the actual convection loss $L=\sum_{i=1}^{T}\parallel \bN(\bv_i)\bv_i-\bN(\tilde{\bv}_i)\bv_i\parallel_{\bM^{-1}}$ is not used in the training session. 

Generally, we observed that the consideration of the FEM norms, e.g., for
computing the loss as $L=\sum_{i=1}^{T}\parallel
\bv_i-\tilde{\bv}_i\parallel_{2}^2$ did not lead to better reconstruction
errors while causing a significant computational overhead. 
Therefore, we sticked to the standard \emph{mean squared error} loss function,
that corresponds to the unweighted $2$-norm of the vectors representing the
states $\bv$ and the forms like $\bN(\bv)$.

% replaced with the $\bM$-norm loss $L=\sum_{i=1}^{T}\parallel
% \bv_i-\tilde{\bv}_i\parallel_{\bM}$, the models show a similar trend of
% reconstruction errors to the results in Figure \ref{fig:rec_err-single} except
% for the CNN model. The CNN is not trained enough for the competition to POD.
% Also, because of the need to include the multiplications by the mass matrix, the
% $M$-norm loss means a significant increase in the training time.

For the iCAE approach, a classification strategy for the velocity states is
learned, once the encoder $\mu$ is available. \Cref{fig:clt-single} shows
that the iCAE5 can capture the periodic flow pattern which is the wake behind a cylinder four seconds later. Literally, data are the most important factor to build data-driven models. In other words, clustering affects reconstruction performance since decoders are trained using clustered datasets. When the states are misclassified, it could cause momentarily sparks that have high errors. 
According to \Cref{fig:v_err-single}, contrary to the overall performance of the iCAE, several abnormal values are observed in the time period $[4,5]$ when $n_\rho=5$. 
%In Figure \ref{fig:traj2}, POD is far behind other models. cPOD5 and iCAE keep up the baseline compared to non-clustering methods. Contrary to the overall performance, several abnormal states are observed in the time period $[3,5]$. The phenomenon could be also estimated in Figure \ref{fig:v_err} and \ref{fig:nvv_err}. It shows that clustering is vulnerable to the range in which the periodic flow pattern starts. 

Next, we use the residual metrics to compare the residuals without % (\ref{metric:residual}, \ref{metric:relresidual})
resorting to the time derivatives.
As mentioned in \Cref{sec:mor}, the clustering models cannot reproduce
$\dot{\tilde{\bv}}$. 
Moreover, as the computation of the data is based on a sophisticated second
order time integration scheme, the error of the numerical evaluation of $\dot v$ only based on
the snapshots that only roughly covers the time domain totally dominated the
residuals.
% to compare fairly the methods. 
% Regarding residuals, the residual of the FOM is not close to zero, because a 2nd order time integration scheme is used.
In \Cref{fig:residual-single}, the errors have a similar trend to the
reconstruction errors shown in \Cref{fig:rec_err-single} with iCAE
performing best and the linear reconstruction by CAE outperforming POD at the
smaller dimensions of $\rho$.
Also not the trajectory plots of the residuals in \Cref{fig:v_err-single} and
\Cref{fig:nvv_err-single}, that clearly show these trends except from the
initial phase where all methods deliver similar approximation results.

%In \Cref{fig:residual-single}, the iCAE follows that of the baseline well. The CAE and the cPOD follow the frequency of the baseline well. 
%However, the values of the CAE are shifted up and the cPOD has a small amplitude in the periodic domain compared to the baseline. 
%The graph of the POD has a proper amplitude and is moved to the right. 
%The residuals of POD-ROM and CAE-ROM

%%%%%%%%%%%%%%%%%%%%%% to rewrite
\subsection{Double Cylinder}\label{sec:num-setup-double-cyl}
%%%%%%%%%%%%%%%%%%%%%%%%%%%%%%%%%%%%%%%%%%%%%%%%%%%%%%%%%%
% Reconstruction Errors
%%%%%%%%%%%%%%%%%%%%%%%%%%%%%%%%%%%%%%%%%%%%%%%%%%%%%%%%%%
\begin{figure}[t]
\centering
\includegraphics[width=1\columnwidth]{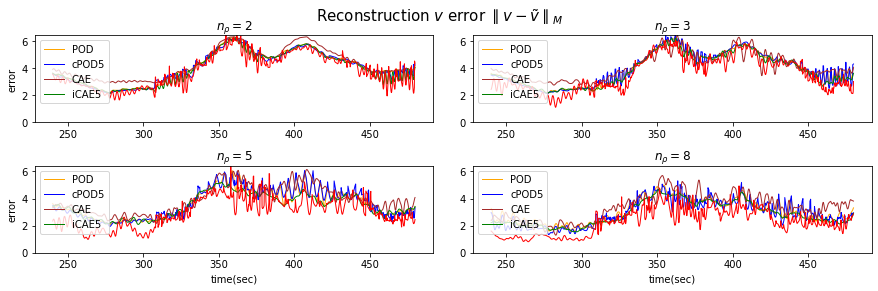}
\caption{Reconstruction error in $[240,480]$ (\Cref{sec:num-setup-double-cyl}).}%
\label{fig:v_err-double}
\end{figure}

\begin{figure}[t]
\centering
\includegraphics[width=1\columnwidth]{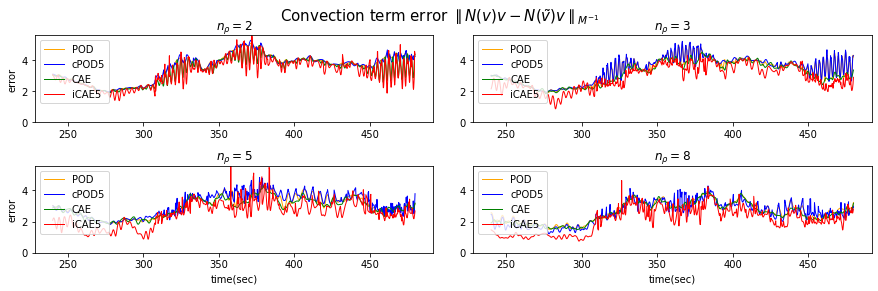}
\caption{Averaged convection error in $[240,480]$ for
the double cylinder case (\Cref{sec:num-setup-double-cyl}).}%
\label{fig:nvv_err-double}
\end{figure}

\begin{figure}[t]
\centering
\includegraphics[width=1\columnwidth]{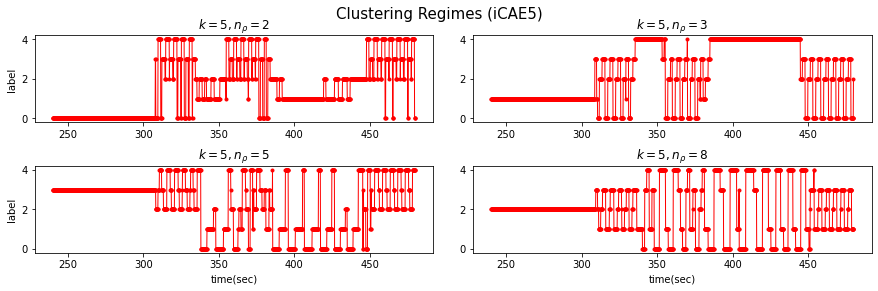}
\caption{Prompt clustering reaction in $[240,480]$ for
the double cylinder case (\Cref{sec:num-setup-double-cyl}).}%
\label{fig:clt-double}
\end{figure}

\begin{figure}[t]
\centering
\includegraphics[width=1\columnwidth]{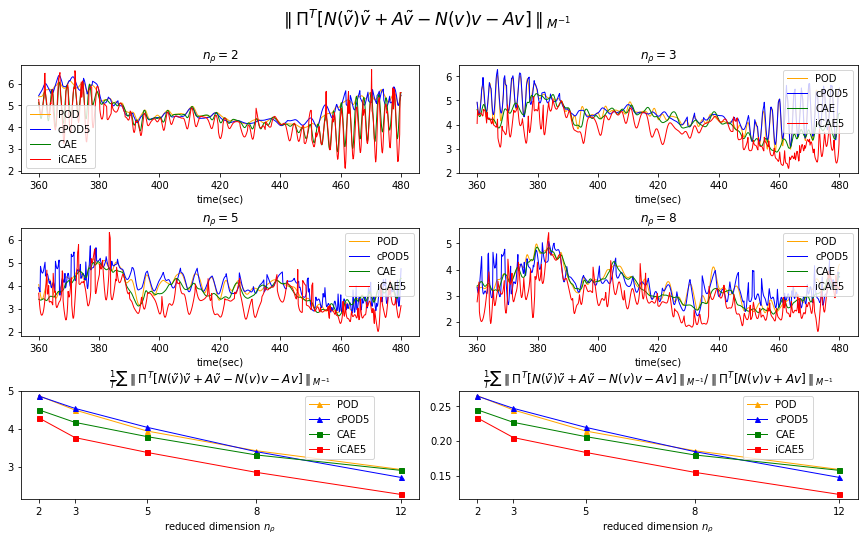}
\caption{Error of $ \Pi^\top(\bN(\bv)\bv+ \bA\bv-\bf)$ in $[360,480]$ (\Cref{sec:num-setup-double-cyl}).}
\label{fig:residual-double}
\end{figure}

As a second example that features richer dynamics, we consider an FEM simulation
of Equation (\ref{ns}) for a double cylinder wake with a Reynolds number of 50;
see \Cref{fig:devsnap} for an illustration of the setup. 
Here, the state vectors $\bv(t)\in\mathbb{R}^{25982}$ (i.e., $n_v=25982$) are
distributed over the physical domain $(-20,70)\times (-20,20)$. Accordingly, to
obtain the CNN input data $\bvcnn(t) \in \mathbb R^{2\times47\times63}$, an
interpolation matrix $\bI_c\in\mathbb{R}^{25982\times (47\cdot 63)}$ is applied. 

For training the models, we use 800 snapshots of $x$ and $y$ directional
velocities in the time interval $[240,480]$, thus leaving aside the startup
phase.

In the evaluation session, 1600 snapshots following the same condition as the training dataset in $[240,480]$ are used. 
The rest of the experimental setting is the same as that of the single cylinder case. As a noncompetitive model, we exclude the CNN model (\ref{cnnnet}) from the experiment and then the existing models with the same hyperparameters are applied for the double-cylinder case.

The averaged errors (see \Cref{fig:rec_err-double}) confirm the findings of the
previous example namely that the iCAE achieve the best results and that the CAE
outperforms POD in particular at low parameter dimensions. 
However, other than for the single cylinder case, a clear gap between the
approaches
in \Cref{fig:v_err-double} and \Cref{fig:nvv_err-double} displaying
the errors over time is not observed. 
This might be due in particular to the many different regimes of the flow that
also appears to be the cause that the clustering results do not show clear
patterns if compared against various dimensions of $\rho$; see
% n contrast with the single-cylinder case, this example contains a chaotic period that can work against clustering that kind of data.
\Cref{fig:clt-double} where the iCAE5 yields different clustering outcomes
depending on $n_\rho=2,3,5,8$. Nonetheless, as shown in
\Cref{fig:residual-double}, the iCAE also achieves the lowest error in the
residuals of $\Pi^\top(\bN(\bv)\bv+ \bA\bv)$ over the reduced dimensions $n_\rho$.

\section{Conclusion}\label{concl}
In this paper, we have investigated several combinations of POD, convolutional
autoencoders, and clustering for low-dimensional parametrizations of fluid flow.
We confirmed that, as suggested by numerous theoretical and numerical studies,
that the nonlinear approaches can outperform linear methods like POD at the very
low dimensions of the parametrization.

Furthermore, we confirmed that neural networks provide a general method for
finding these low-dimensional parametrization even in the high-dimensional data
regime, in particular if one can exploit
sparsity patterns as they are implicitly contained in convolutional neural
networks. 

In order to probe the limits of the parametrization, we have included $k$-means
clustering on the reduced coordinates so that the reconstruction can be made
specific to specific flow regimes. 
Although this approach (in terms of the iCAE method
in this paper) reaches significantly lower error levels in particular at low
dimensions, its practical use, e.g., in simulations, is in question because of
the nonsmooth operations for the decoding.

In the future, we will investigate the performance of the proposed
low-dimensional parametrizations in simulations and controller design. 
Another relevant future development will concern relaxations of the clustering
maps so that reconstruction based on clustering information will become a smooth
and possibly even linear operation.

\section{Code Availability}\label{sec:code-availability}

%Upon request, we can provide the complete code base of the presented approaches
%and reported numerical results to the reviewers. Upon acceptance, we will make
%the code publicly available.

The source code of the implementations used to compute the presented results is available from 
  \href{https://doi.org/10.5281/zenodo.7575808}{\texttt{doi:10.5281/zenodo.7575808}}
under the Creative Commons Attribution 4.0 international license and is authored
by Yongho Kim.

\section*{Acknowledgments}%
\addcontentsline{toc}{section}{Acknowledgments}
We acknowledge funding by the German Research Foundation (DFG) through the research training group 2297 ``MathCoRe'', Magdeburg.

\newpage
\bibliographystyle{abbrvurl}
\bibliography{cae-cluster-pod}
% \bibliography{cae-cluster-pod,csc}

%\section{Appendix}

\end{document}